\documentclass[12 pt]{amsart}
\usepackage{amscd,amssymb,amsmath,amsthm}

\newtheorem{Theorem}{Theorem}

\newcommand{\rarr}{\rightarrow}

\newcommand{\com}{\mathbb{C}}

\newcommand{\M}{{\mathcal{M}}}

\begin{document}
\baselineskip=16pt

\begin{center}
{{\bf The $\kappa$ classes on the moduli spaces of curves}}
\end{center}

\vspace{15pt}

\begin{center}
{{\em R. Pandharipande}}
\end{center}

\vspace{15pt}

\begin{center}
{{\em August 2011}}
\end{center}

$$ $$

In the past few years, substantial
progress has been made in the understanding
of the algebra of  $\kappa$ classes on the moduli spaces of curves.
My goal here is to provide a short introduction to
the new results. Along the way, I will discuss several open questions.
The article accompanies my talk at {\em A celebration of algebraic
geometry} at Harvard
in honor of the $60^{th}$ birthday of J. Harris.

\vspace{15pt}
\noindent {\bf A. Moduli spaces of curves}
\vspace{15pt}

Let $\M_g$ be the moduli space of compact nonsingular 
curves over $\com$
of genus $g$. We view $\M_g$ as a nonsingular Deligne-Mumford
stack of dimension $3g-3$. We will consider the spaces
\begin{equation}\label{geet}
\M_g \subset \M_g^c \subset \overline{\M}_g\ . 
\end{equation}
Here, $\M_g^c$ is the moduli space of curves of compact type 
(curves with no cycles in the dual graph). While $\M_g$ and
$\M_g^c$ are open, the moduli space of stable curves
$\overline{\M}_g$ is compact.
The Chow rings of the moduli spaces \eqref{geet} are well-defined. 
Because of stack considerations, we will take the
Chow rings with $\mathbb{Q}$-coefficients.

For the investigation of the $\kappa$ classes, the moduli
spaces of curves with markings 
\begin{equation*}
\M^{rt}_{g,n} \subset \M_{g,n}^c \subset \overline{\M}_{g,n}\ . 
\end{equation*}
play an essential role and should be treated on the same footing.
 The moduli space 
$\M^{rt}_{g,n}$ of curves with rational tails is the inverse
image of $\M_g$ under the forgetful map
$$\overline{\M}_{g,n} \rightarrow \overline{\M}_{g}\ . $$

\pagebreak
\vspace{15pt}
\noindent {\bf B. $\kappa$ classes}
\vspace{15pt}

The $\kappa$ classes in the Chow ring 
$A^*(\overline{\M}_{g,n})$
are defined by the following geometry.
Let
$$\epsilon: \overline{\M}_{g,n+1} \rarr \overline{\M}_{g,n}$$
be the universal curve viewed as the ($n+1$)-pointed space, let
$$\mathbb{L}_{n+1} \rarr \overline{\M}_{g,n+1}$$
be the line bundle obtained from the cotangent space
of the last marking, 
and let
$$\psi_{n+1} = c_1(\mathbb{L}_{n+1})\ \in A^1(\overline{\M}_{g,n+1})$$
be the Chern class.
The $\kappa$ classes, first defined by Mumford, 
are
$$\kappa_i = \epsilon_*(\psi^{i+1}_{n+1})\  \in A^i(\overline{\M}_{g,n}), 
\ \ \ i \geq 0\ .$$
The simplest is $\kappa_0$ which equals $2g-2+n$ times the unit
in $A^0(\overline{\M}_{g,n})$.
The convention 
$$\kappa_{-1}= \epsilon_*(\psi_{n+1}^0)=0$$
is often convenient. 

The $\kappa$ classes on $\M^{rt}_{g,n}$ and $\M_{g,n}^c$
are defined via restriction from $\overline{\M}_{g,n}$.
Define the $\kappa$ rings
\begin{eqnarray*}
\kappa^*(\M^{rt}_{g,n}) & \subset & A^*(\M^{rt}_{g,n}),\\
\kappa^*(\M_{g,n}^c)& \subset &  A^*(\M_{g,n}^c), \\
\kappa^*(\overline{\M}_{g,n}) & \subset & A^*(\overline{\M}_{g,n}),
\end{eqnarray*}
to be the $\mathbb{Q}$-subalgebras 
generated by the $\kappa$ classes.
Of course,  the $\kappa$ rings are graded by degree.

Since $\kappa_i$ is a  tautological class
the $\kappa$ rings
are subalgebras of the corresponding tautological rings.
For unpointed nonsingular curves, the $\kappa$ ring equals
the tautological ring by definition \cite{Mum},
$$\kappa^*(\M_g)= R^*(\M_g)\ .$$
Otherwise, the inclusion of the $\kappa$ ring in the
tautological ring is usually proper.

\pagebreak

\vspace{15pt}
\noindent {\bf C. $M_g$ and the Faber-Zagier conjecture}
\vspace{15pt}

Consider first the ring $\kappa^*(\M_g)$. The basic non-vanishing
and vanishing results,
$$\kappa^{g-2}(\M_g) \cong \mathbb{Q}, \ \ \
\kappa^{>g-2}(\M_g)=0\ $$
have been known now for some time \cite{Faber, L}. 
Moreover, there are no relations \cite{Bold} among the $\kappa$ classes
of degree less than or equal to ${\lfloor\frac{g}{3}\rfloor}$, and the
classes $$\kappa_1, \ldots, \kappa_{{\lfloor\frac{g}{3}\rfloor}}$$ generate
$\kappa^*(\M_g)$ \cite{Ionel, Morita}. 
These properties had all been conjectured earlier by Faber \cite{Faber}.

We define a
set of relations as follows.
Let 
$$\mathbf{p} = \{\ p_1,p_3,p_4,p_6,p_7,p_9,p_{10}, \ldots\ \}$$
be a variable set indexed by positive integers not congruent
to $2$ mod 3.
Let
\begin{multline*}
\Psi(t,\mathbf{p}) =
(1+tp_3+t^2p_6+t^3p_9+\ldots) \sum_{i=0}^\infty \frac{(6i)!}{(3i)!(2i)!} t^i
\\ +(p_1+tp_4+t^2p_7+\ldots) 
\sum_{i=0}^\infty \frac{(6i)!}{(3i)!(2i)!} \frac{6i+1}{6i-1} t^i \ .
\end{multline*}
Define the constants $C^r(\sigma)$ by the formula
$$\log(\Psi)= 
\sum_{\sigma}
\sum_{r=0}^\infty C^r(\sigma)\ t^r 
\mathbf{p}^\sigma
\ . $$
Here, $\sigma$ denotes a partition of size 
$|\sigma|$ which avoids all 
 parts congruent to 2 mod 3. If $\sigma=1^{n_1}3^{n_3}4^{n_4} \ldots$,
then $\mathbf{p}^\sigma= p_1^{n_1}p_3^{n_3}p_4^{n_4}\ldots$ as usual.
Let 
$$\gamma= 
\sum_{\sigma}
 \sum_{r=0}^\infty C^r(\sigma)
\ \kappa_r t^r 
\mathbf{p}^\sigma
\ .
$$

\begin{Theorem} 
In $\kappa^r(\M_g)$, the relation
$$
\big[ \exp(-\gamma) \big]_{t^r \mathbf{p}^\sigma}  = 0$$
holds when
$g-1+|\sigma|< 3r$ and
$g\equiv r+|\sigma|+1 \mod 2$.
\end{Theorem}

The relations of Theorem 1, called the {\em FZ relations},
 were conjectured to hold several
years ago by Faber and Zagier
from low genus data and a study of the
Gorenstein quotient of $\kappa^*(\M_g)$. Guessing the full
structure here was certainly a remarkable feat.
Theorem 1 was proven by myself and A. Pixton \cite{PP} last
year using the geometry of stable quotients \cite{MOP}.

To the best of our knowledge, a relation in $\kappa^*(\M_g)$ which is
not in the span of the FZ relations has not yet been found.
In particular, all relations obtained from the various geometrical
constructions attempted in the past \cite{Faber} appear to be
covered by Theorem 1.
Whether Theorem 1 exhausts all relations in $\kappa^*(\M_g)$ is
a very interesting question.

\vspace{+10pt}
\noindent{\bf Q1.} Are all relations among the $\kappa$
classes in $\kappa^*(\M_g)$ generated by\\
\hspace*{22pt}    Theorem 1?
\vspace{+10pt}

Theorem 1 only provides finitely many relations
in $\kappa^r(\M_g)$ for fixed $g$ and $r$, and thus may be 
calculated completely. 
When the relations yield a Gorenstein
ring with socle in $\kappa^{g-2}(\M_g)$, no further relations are possible.
However, the relations of Theorem 1 do not always yield such a
Gorenstein ring (failing first in genus 24 as checked by Faber).
For $g<24$, Faber's calculations show
Theorem 1 does provide all relations in $\kappa^*(\M_g)$.
For higher genus $g\geq 24$, either Theorem 1 fails to
provide all the relations in $\kappa^*(\M_g)$ {or}
$\kappa^*(\M_g)$ is not Gorenstein. 

\vspace{+10pt}
\noindent{\bf Q2.} Is $\kappa^*(M_g)$ a Gorenstein ring
with socle in degree $g-2$?
\vspace{+10pt}

Faber's original conjecture \cite{Faber} asserts an affirmative
answer to Q2. Questions Q1 and Q2 can not both have an affirmative
answer in genus 24. Which assertion is false?

\vspace{10pt}
The main actors in the FZ relations are the functions 
\begin{eqnarray*}
A(z)& = &\sum_{i=0}^\infty \frac{(6i)!}{(3i)!(2i)!} 
\left(\frac{z}{72}\right)^i\ , \\
B(z)& = & \sum_{i=0}^\infty \frac{(6i)!}{(3i)!(2i)!} \frac{6i+1}{6i-1} 
\left(\frac{z}{72}\right)^i \ 
\end{eqnarray*}
written here in the variable $z=72t$.
The function $B$ is determined from $A$ by the 
differential equation
\begin{equation*}
-\frac{1}{2}A + zA + 6 z^2\ \frac{dA}{dz} = \frac{1}{2}B\ .
\end{equation*}
The main hypergeometric differential equation satisfied
by $A$ is
\begin{equation}\label{r4r}
36 z^2 \frac{d^2}{dz^2} A + (72z-6) \frac{d}{dz} A + 5 A = 0 \ .
\end{equation}
A more open ended question is following.

\vspace{+10pt}
\noindent{\bf Q3.} What is the meaning of the function $A$? Does $A$
or the\\ \hspace*{22pt} differential equation \eqref{r4r} occur elsewhere
     in mathematics?
\vspace{+10pt}

The statement of Theorem 1 contains a peculiar mod 2 condition.
As a consequence, the relations in $\kappa^*(\M_g)$ given by Theorem 1
also hold in $\kappa^*(\M_{g-2})$. Note $\kappa_0$ specializes to
$2g-2$ in genus $g$ and $2g-6$ in genus $g-2$ !  
In our proof \cite{PP} of Theorem 1, the mod 2 condition
arises (after a considerable amount
of work) from the 2 fixed points of $\com^*$ on $\mathbb{P}^1$.

\vspace{+10pt}
\noindent{\bf Q4.} Is there a simple explanation of the mod 2 condition?
\vspace{+10pt}

We can also ask about relations on the moduli 
space $\M_{g,n}^{rt}$ with markings.
A. Pixton has a precise proposal for the $n=1$ case in the form
of Theorem 1. Work here is just starting and will be reported by Pixton 
elsewhere. 
 
\vfill

\pagebreak

\vspace{15pt}
\noindent {\bf D. Universality for $M_{g,n}^c$}
\vspace{15pt}

Consider next the ring $\kappa^*(\M_{g,n}^c)$. The basic non-vanishing
and vanishing results,
$$\kappa^{2g-3+n}(\M^c_{g,n}) \cong \mathbb{Q}, \ \ \
\kappa^{>2g-3+n}(\M^c_{g,n})=0\ $$
have been proven in \cite{FPrel,GV}. The classes
$$\kappa_1,\kappa_2,\ldots, \kappa_{g-1+\lfloor\frac{n}{2}
\rfloor}$$
generate $\kappa^*(\M_{g,n}^c)$.
If $n>0$, there are no relations of degree less than or equal to
${g-1+\lfloor\frac{n}{2}
\rfloor}$. In $\kappa^*(\M_g^c)$, there are no relations of degree
less than $g-1$ (whether degree $g-1$ relations can occur is
not known).
The proofs of the above generation and freeness results
can be found in \cite{kap}.

A surprising feature about the $\kappa$ rings in the compact
type case is the following universality result proven in \cite{kap}.

\begin{Theorem}
Let $g>0$ and $n>0$, then the assignment $\kappa_i \mapsto \kappa_i$
extends to a ring isomorphism
$$ \iota: \kappa^*(M_{g-1,n+2}^{c}) \cong  \kappa^*(M_{g,n}^c) \ . $$
\end{Theorem}

In other words, the relations among the $\kappa$ classes
in the above cases 
are genus independent. 
By composing the isomorphisms $\iota$ of Theorem 2, we obtain isomorphisms
$$\iota: \kappa^*(M_{0,2g+n}^c) \cong \kappa^*(M_{g,n}^c)$$
so long as $n>0$. Hence, universality reduces all questions about
the $\kappa$ rings to genus 0. 

In \cite{kap}, calculations of the relations,
bases, and Betti numbers of the ring
$\kappa^*(M_{g,n>0}^c)$ are 
obtained using the genus 0 reduction.
Let $P(d)$ be the set of partitions of $d$, and let
$$P(d,k)\subset P(d)$$ be the set of partitions of $d$ into
at most $k$ parts. 
Let $|P(d,k)|$ be the
cardinality. To a partition with positive parts
$\mathbf{p}= (p_1,\ldots,p_\ell)$ in $P(d,k)$,
we associate a $\kappa$ monomial by
$$\kappa_{\mathbf{p}} = \kappa_{p_1} \cdots \kappa_{p_\ell} \in 
\kappa^d(M_{g,n}^c) \ .$$

\begin{Theorem} \label{htt5} For $n>0$,
a $\mathbb{Q}$-basis of $\kappa^d(M_{g,n}^c)$ is given by
$$\{ \kappa_{\mathbf{p}} \ | \ \mathbf{p} \in P(d,2g-2+n-d)\ \} \  .$$
\end{Theorem}

The Betti number calculation,
$$\text{dim}_{\mathbb{Q}} \ \kappa^d(M_{g,n}^c) \ = \ |P(d,2g-2+n-d)|\ ,$$
is of course implied by Theorem \ref{htt5}.

\vspace{+10pt}
\noindent{\bf Q5.} Is there an analogue of Schubert calculus
in the $\kappa$ ring 
with respect \\
\hspace*{22pt}  to the basis of Theorem \ref{htt5}?
\vspace{+10pt}

The tautological rings $R^*(M^c_{g,n})$ have been
conjectured in \cite{FP1,Pand} to be Gorenstein algebras with socle in degree
$2g-3+n$,
$$\phi: R^{2g-3+n}(M_{g,n}^c) \stackrel{\sim}{\rarr} \mathbb{Q}\ .$$
The following result \cite{kap}
shows the socle evaluation is as non-trivial
as possible on the $\kappa$ ring. 

\begin{Theorem} \label{nmmg}
If $n>0$ and $\xi \in \kappa^d(M_{g,n}^c) \neq 0$, the linear
function 
$$L_\xi: R^{2g-3+n-d}(M_{g,n}^c)\rarr \mathbb{Q}$$
defined by the socle evaluation
$$L_\xi(\gamma) =\phi(\gamma\cdot \xi)$$
is non-trivial.
\end{Theorem}

All of the above results for the $\kappa$ rings in the compact type
case require at least 1 marked point. In the unpointed case $n=0$,
half of the universality still holds: the
assignment $\kappa_i\mapsto \kappa_i$ extends to a surjection
$$\iota_g: \kappa^*(M_{0,2g}^c) \rightarrow \kappa^*(M_{g}^c)\ . $$
However, a non-trivial kernel is possible. The first kernel
occurs in genus $g=5$.

\vspace{+10pt}
\noindent{\bf Q6.} What is the kernel of
$\kappa^*(M_{0,2g}^c) \rightarrow \kappa^*(M_{g}^c)$?
\vspace{+10pt}

In genus 5, the kernel of $\iota_5$ is related to Getzler's
relation in $\overline{M}_{1,4}$, see \cite{kap} for a discussion.
A complete answer to Q6  will likely involve sequences
of special relations in the
tautological ring.

\vspace{10pt}
A natural question to ask is whether the $\kappa$ relations
in the compact type case can be put in a  form parallel to Theorem 1.
An affirmative answer has been found by A. Pixton.

We define a
set of relations as follows.
Let 
$$\mathbf{p} = \{\ p_1,p_2,p_3,\ldots\ \}$$
be a variable set indexed by all positive integers.
Let
\begin{multline*}
\Psi(t,\mathbf{p}) =
(1+tp_2+t^2p_4+t^3p_6+\ldots) \sum_{i=0}^\infty (2i-1)!!\ t^i
\\ +(p_1+tp_3+t^2p_5+\ldots) \  
\end{multline*}
where $(2i-1)!!= \frac{(2i)!}{2^ii!}$ as usual.
Define the constants $C^r(\sigma)$ by the formula
$$\log(\Psi)= 
\sum_{\sigma}
\sum_{r=0}^\infty C^r(\sigma)\ t^r 
\mathbf{p}^\sigma
\ . $$
Here, $\sigma$ denotes any partition with positive parts.
Let 
$$\gamma= 
\sum_{\sigma}
 \sum_{r=0}^\infty C^r(\sigma)
\ \kappa_r t^r 
\mathbf{p}^\sigma
\ .
$$

\begin{Theorem} {\em [Pixton]}
In $\kappa^r(\M_{g,n}^c)$, the relation
$$
\big[ \exp(-\gamma) \big]_{t^r \mathbf{p}^\sigma}  = 0$$
holds when
$2g-2+n+|\sigma|<2r$.
\end{Theorem}

Pixton's proof uses the Gorenstein criterion of Theorem 4
to check the relations. Furthermore, Pixton proves the
relations of Theorem 5 generate {\em all} the relations
among $\kappa$ classes in $\kappa^*(\M_{g,n}^c)$ 
in case $n>0$.

\vspace{15pt}
\noindent {\bf E. Moduli of stable curves}
\vspace{15pt}

The socle evaluations of polynomials in the $\kappa$ classes
are known for all three cases
$$\kappa^*(\M_{g,n}^{rt}), \ \ \kappa^*(\M_{g,n}^c), \ \
\kappa^*(\overline{M}_{g,n})\ . $$ 
The $\kappa$ evaluations can be transformed to 
Hodge integrals. The relevant evaluations of the
latter are summarized in \cite{FP1}.
The socle evaluation in the stable curve case is
obtained from the proof of Witten's conjecture \cite{Wit}.

The socle evaluations imply relations in the top degree
$$\kappa^{g-2+n}(\M_{g,n}^{rt}), \ \ \kappa^{2g-3+n}(\M_{g,n}^c), \ \
\kappa^{3g-3+n}(\overline{M}_{g,n})\ . $$ 
The investigation of the $\kappa$ ring concerns
relations in all degrees. In the stable curve case,
no uniform results are known at the moment for 
$\kappa$ relations above the socle.

\vspace{+10pt}
\noindent{\bf Q7.}
Is there a formula for relations in $\kappa^*(\overline{\M}_{g,n})$
parallel to \\
\hspace*{22pt}    
Theorem 1 and Theorem 5?
\vspace{+10pt}

\vspace{15pt}
\noindent {\bf Acknowledgments}
\vspace{15pt}

Most of my work on $\kappa$ classes
is related to conversations over
the years with C. Faber whom I met for the first time on
a spring morning in Joe's office when I was a graduate
student. Many of the results reported here were
found in collaboration with
A. Pixton (or by Pixton himself)
in Princeton and in Lisbon last year. My paper with A. Marian and D. Oprea
on the moduli space of stable quotients plays an
essential role in the proofs.
At the Instituto
Superior T\'ecnico in Lisbon, I was supported by a Marie Curie fellowship and
a grant from the Gulbenkian foundation. 
I am supported by NSF grant DMS-1001154.

\vspace{+20 pt}
\noindent Departement Mathematik \hfill Department of Mathematics \\
\noindent ETH Z\"urich \hfill  Princeton University \\
\noindent 8092 Z\"urich  \hfill Princeton, NJ 08544 \\
\noindent Switzerland  \hfill USA \\


\begin{thebibliography}{[(III)]}


\bibitem {Bold} S. Boldsen, {\em Improved homological stability for the
mapping class group with integral or twisted coefficients},
arXiv:0904.3269.


\bibitem {Faber}
C.~Faber, 
\newblock{\em A conjectural description of the tautological ring of the moduli s
pace of curves}, Moduli of curves and abelian varieties,  109--129, Aspects Math., 
Vieweg, Braunschweig, 1999.


\bibitem{FP1} C. Faber and R. Pandharipande (with an appendix
by D. Zagier), {\em Logarithmic series
and Hodge integrals in the tautological ring}, Michigan Math. J. {\bf 48}
(2000), 215--252.


\bibitem{FPrel}
C. Faber and R. Pandharipande, {\em Relative maps and tautological
classes}, JEMS {\bf 7} (2005), 13--49.

\bibitem {GV} T. Graber and R. Vakil, {\em                                      
Relative virtual localization and vanishing of tautological classes on moduli s\
paces of curves}, Duke Math. J. {\bf 130} (2005), 1–37.


\bibitem{Ionel} E. Ionel, {\em Relations in the tautological ring
of $\M_g$}, Duke Math. J. {\bf 129} (2005), 157--186.

\bibitem{L} E. Looijenga, {\em
On the tautological ring of $M_g$}.  
Invent. Math. {\bf 121}  (1995), 411--419.


\bibitem{MOP} 
A.~Marian, D. ~Oprea, and R.~Pandharipande, \newblock{\em 
The moduli space of stable quotients}, Geom. \& Top. (to appear).



\bibitem {Morita} S. Morita, {\em Generators for the
tautological algebra of the moduli space of curves},
Topology {\bf 42} (2003), 787--819.

\bibitem{Mum} D. Mumford, {\em Towards an enumerative geometry of
the moduli space of curves}, in {\em Arithmetic and Geometry}
(M. Artin and J. Tate, eds.), Part II, Birkh\"auser, 1983, 271-328.

\bibitem{Pand}
R. Pandharipande, {\em Three questions in Gromov-Witten theory},
Proceedings of the ICM (Beijing 2002), Vol. II, 503--512.


\bibitem{kap} R. Pandharipande, {\em The $\kappa$ ring of the
moduli of curves of compact type}, Acta Math. (to appear).


\bibitem{PP} R. Pandharipande and A. Pixton, {\em Relations
in the tautological ring}, arXiv:1101.2236.

\bibitem{Wit} E. Witten, {\em Two dimensional gravity and intersection
theory on moduli space}, Surveys in Diff. Geom. {\bf 1} (1991), 243--310.



\end{thebibliography}
\end{document}